# DISCUSSION: THE DANTZIG SELECTOR: STATISTICAL ESTIMATION WHEN $p$ IS MUCH LARGER THAN $n$

By T. Tony Cai[1] and Jinchi Lv

*University of Pennsylvania and Princeton University*

Professors Candès and Tao are to be congratulated for their innovative and valuable contribution to high-dimensional sparse recovery and model selection. The analysis of vast data sets now commonly arising in scientific investigations poses many statistical challenges not present in smaller scale studies. Many of these data sets exhibit sparsity where most of the data corresponds to noise and only a small fraction is of interest. The needs of this research have excited much interest in the statistical community. In particular, high-dimensional model selection has attracted much recent attention and has become a central topic in statistics. The main difficulty of such a problem comes from collinearity between the predictor variables. It is clear from the geometric point of view that the collinearity increases as the dimensionality grows.

A common approach taken in the statistics literature is the penalized likelihood, for example, Lasso (Tibshirani [11]) and adaptive Lasso (Zou [12]), SCAD (Fan and Li [7] and Fan and Peng [9]) and nonnegative garrote (Breiman [1]). Commonly used algorithms include LARS (Efron, Hastie, Johnstone and Tibshirani [6]), LQA (Fan and Li [7]) and MM (Hunter and Li [10]). In the present paper, Candès and Tao take a new approach, called the Dantzig selector, which uses $\ell_1$-minimization with regularization on the residuals. One promising fact is that the Dantzig selector solves a linear program, usually faster than the existing methods. In addition, the authors establish that, under the Uniform Uncertainty Principle (UUP), with large probability the Dantzig selector mimics the risk of the oracle estimator up to a logarithmic factor $\log p$, where $p$ denotes the number of variables.

We appreciate the opportunity to comment on several aspects of this article. Our discussion here will focus on four issues: (1) connection to sparse signal recovery in the noiseless case; (2) the UUP condition and identifiability of the model; (3) computation and model selection; (4) minimax rate.

Received January 2007.
[1]Supported in part by NSF Grant DMS-06-04954.

**1. Sparse signal recovery.** The "large $p$, small $n$" regression problem considered in this paper can be viewed as a generalization of the classical linear algebra problem in which one wishes to solve the linear equation

$$y = X\beta, \tag{1}$$

where $X$ is a given $n \times p$ matrix and $y$ is a vector in $\mathbb{R}^n$. Because $p > n$, the linear equation (1) is underdetermined and there are an infinite number of solutions to the equation. The goal is to find the "sparsest" solution under certain regularity conditions. In this noiseless setting, the Dantzig selector reduces to an $\ell_1$-minimization over the space of all representations of the signal:

$$\text{Minimize } \|\beta\|_1 \quad \text{subject to} \quad y = X\beta. \tag{2}$$

This idea of finding a sparse representation using $\ell_1$-minimization has been used in Donoho and Elad [4]. The authors have also used this approach in their earlier work (Candès and Tao [2, 3]) on recovering sparse signals in the noiseless case. When adding a Gaussian noise term $\varepsilon$ to (1), the linear algebra problem becomes a nonstandard linear regression problem because $p \gg n$. In the classical linear regression problem when $p \leq n$ the least squares estimator is the solution to the normal equation

$$X^T y = X^T X \beta. \tag{3}$$

The constraint $\|X^T r\|_\infty \leq \lambda_p \sigma$ (which is the same as $\|X^T y - X^T X \beta\|_\infty \leq \lambda_p \sigma$) in the convex program (DS) in this paper can be viewed as a relaxation of the normal equation (3). And similarly to the noiseless case $\ell_1$-minimization leads to the "sparsest" solution over the space of all feasible solutions.

The authors suggest using $\lambda_p = \sqrt{2 \log p}$, which is equal to $\sqrt{2 \log n}$ in the orthogonal design setting. In this setting, the oracle properties of the Dantzig selector are in line with those shrinkage results in Donoho and Johnstone [5] which are shown to be optimal in the minimax sense. When $p \gg n$, it might be possible that the regularization factor $\lambda_p = \sqrt{2 \log p}$ in the Dantzig selector overshrinks the $p$-vector $\beta$ and underestimates the nonzero coordinates. It would be interesting to find the "optimal" regularization factor.

**2. The UUP condition and identifiability of the model.** A nice idea in this paper is the use of the UUP condition. The UUP condition has also been used by the authors in their earlier work (Candès and Tao [2, 3]) in the noiseless setting. The UUP condition roughly says that for any small set of predictors, these $n$-vectors are nearly orthogonal to each other. The authors give an interpretation of the UUP condition in terms of the model identifiability and have established oracle inequalities for the Dantzig selector under the UUP using geometric arguments. However, we still have some concerns about this condition.

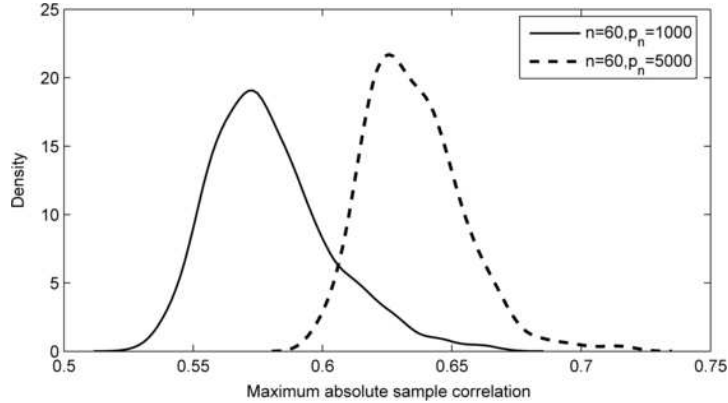

Fig. 1. *Distribution of the maximum absolute sample correlation when $n = 60$, $p_n = 1000$ (solid curve) and $n = 60$, $p_n = 5000$ (dashed curve).*

First, it is computationally unrealistic to verify whether the UUP condition holds for a given design matrix $X$ when $p$ is large and the number of signals $S$ is not too small. Note that computing the $S$-restricted isometry constant $\delta_S$ of the design matrix $X$ is over the space of all $S$-subsets of $\{1, \ldots, p\}$, which is of cardinality $\binom{p}{S}$. This combinatorial complexity makes it infeasible to check the UUP condition for reasonable values of $p$ and $S$, say $p = 1000$ and $S = 5$. So it is interesting to look for other checkable conditions that are compatible with the model identifiability condition in a similar or new setup.

Second, as discussed in Fan and Lv [8], the UUP is hard to satisfy when the dimension $p$ grows rapidly with the sample size $n$. This is essentially due to significant sample correlation, that is, strong collinearity, between the predictor variables in the high-dimensional setting of $p \gg n$. For instance, we take $p$ independent predictors $X_1, \ldots, X_p$ from the standard Gaussian distribution and compute the maximum of the pairwise absolute sample correlations from an $n \times p$ design matrix $X$. Figure 1, which is extracted from Fan and Lv [8], shows the distributions of the maximum correlation with $n = 60$, $p = 1000$ and $n = 60$, $p = 5000$, respectively. The maximum sample correlation between predictors can be very large and close to 1. Moreover, the maximum of the first canonical correlations between two groups of predictors, for example, three predictors in one group and five in another, can be much larger since there are $\binom{p}{3}\binom{p-3}{5} = O(p^8)$ ($\gg n$) choices in this example.

**3. Computation and model selection.** Due to its nature of involving linear programming, the Dantzig selector can be solved quickly and efficiently by a primal–dual interior point algorithm when the dimension is not ultrahigh, for example, in the thousands. It is usually faster to implement than other existing methods such as Lasso. However, in problems of large or ultralarge scale the computational cost of implementing linear programs is still a potential hurdle. For example, in the analysis of microarray gene expression or proteomics data, it is common to have dimension of tens of thousands and implementing linear programs in such settings can still be computationally challenging. Therefore, it is interesting and necessary to study ultrahigh-dimensional model selection. Recently, Fan and Lv [8] introduced a procedure for screening variables via iteratively thresholded ridge regression and proposed a new method of dimension reduction called Sure Independence Screening, for ultrahigh-dimensional feature space. The method can improve the estimation accuracy significantly while speeding up variable selection drastically.

In our experience, we have found that the algorithm solving the Dantzig selector is sensitive to the initial value. Trivial initial values such as constant vectors usually do not work well. The generalized least squares estimator $(X^T X)^- X^T y$ can be used as an initial value. However, the solution is usually nonsparse in our experience.

As mentioned earlier, the regularization factor $\lambda_p = \sqrt{2 \log p}$ in the Dantzig selector leads to relatively large bias in estimating the sparse regression coefficient vector $\beta$. To reduce the bias, the authors suggest a two-stage procedure called the Gauss–Dantzig selector which uses the original Dantzig selector for variable selection and then runs ordinary least squares on the selected variables. It would be interesting to know the theoretical and numerical properties of this and other variations.

**4. Minimax rate.** It is appealing that the Dantzig selector achieves within a logarithmic factor $\log p$ of the ideal risk. We are curious about the optimality of this factor. It is unclear at this point whether the minimax factor should be $\log p$ or $\log n$ or some other quantity. When $p$ is polynomial in the sample size $n$, $\log p$ and $\log n$ are of the same order and the difference between the two is not significant. However, when $p$ is exponential in $n$, say $p = e^{n^a}$ for some $a > 0$, then $\log p = n^a$ becomes large, much larger than $\log n$. It is of theoretical and practical interest to study the minimax behavior of the problem. If the minimax rate is $\log p$, then the Dantzig selector is a rate optimal minimax procedure. Otherwise, it is interesting to construct a procedure that can attain the minimax rate.

**5. Concluding remarks.** $\ell_1$-regularization in terms of linear programming provides a new framework for model selection and is proven effective in solving high-dimensional problems. The Dantzig selector provides us new insight on high-dimensional model selection. Clearly, there is much work ahead of us. This paper opens a new area and will definitely stimulate new ideas in the future. We thank the authors for their clear and imaginative work.

DEPARTMENT OF STATISTICS
THE WHARTON SCHOOL
UNIVERSITY OF PENNSYLVANIA
PHILADELPHIA, PENNSYLVANIA 19104
USA
E-MAIL: tcai@wharton.upenn.edu

DEPARTMENT OF MATHEMATICS
PRINCETON UNIVERSITY
PRINCETON, NEW JERSEY 08544
USA
E-MAIL: jlv@princeton.edu